\documentclass{elsart}
\usepackage{graphicx}
\usepackage{amssymb}
\usepackage{amsmath}
\usepackage{amsfonts}
\begin{document}
\begin{frontmatter}
\title{BDDC by a Frontal Solver and
the Stress Computation in a Hip Joint Replacement}
\author[a,d]{Jakub~\v{S}\'{\i}stek},
\ead{jakub.sistek@fs.cvut.cz}
\author[b]{Jaroslav~Novotn\'{y}},
\ead{novotny@it.cas.cz}
\author[c]{Jan Mandel},
\ead{jan.mandel@cudenver.edu}
\author[a]{Marta \v{C}ert\'{\i}kov\'{a}},
\ead{marta.certikova@fs.cvut.cz}
\author[a]{Pavel Burda}
\ead{pavel.burda@fs.cvut.cz}
\address[a]{Department of Mathematics,
Faculty of Mechanical Engineering \\ Czech Technical University in Prague}
\address[b]{Department of Mathematics,
Faculty of Civil Engineering \\ Czech Technical University in Prague}
\address[c]{Department~of~Mathematical~and~Statistical~Sciences,
University~of~Colorado~Denver}
\address[d]{Institute~of~Thermomechanics,
Academy~of~Sciences of the~Czech~Republic}
\begin{abstract}
A parallel implementation of the BDDC method using the frontal solver
is employed to solve systems of linear equations from finite element
analysis, and incorporated into a standard finite element system for
engineering analysis by linear elasticity. Results of computation of
stress in a hip replacement are presented. The part is made of
titanium and loaded by the weight of human body.
The performance of BDDC with added constraints by averages and with added corners is compared.
\end{abstract}
\end{frontmatter}

\section{Introduction}

\label{sec:Introduction}

Parallel numerical solution of linear problems arising from linearized
isotropic elasticity discretized by finite elements is important in many areas
of engineering. The matrix of the system is typically large, sparse, and
ill-conditioned. The classical frontal solver \cite{Irons-1970-FSS} has became
a popular direct method for solving problems with such matrices arising from
finite element analyses. However, for large problems, the computational cost
of direct solvers makes them less competitive compared to iterative methods,
such as the preconditioned conjugate gradients (PCG) \cite{Concus-1976-GCG}.
The goal is then to
design efficient preconditioners that result in a lower overall cost and can
be implemented in parallel, which has given rise to the field of domain
decomposition and iterative substructuring \cite{Toselli-2005-DDM}.

The Balancing Domain Decomposition based on Constraints (BDDC)
\cite{Dohrmann-2003-PSC} is one of the most advanced preconditioners of this
class. However, the additional custom coding effort required can be an
obstacle to the use of the method in an existing finite element code. We
propose an implementation of BDDC built on top of common components of
existing finite element codes, namely the frontal solver and the element
stiffness matrix generation. The implementation requires only a minimal amount
of additional code and it is therefore of interest. For an important
alternative implementation of BDDC, see~\cite{Li-2006-FBB}.

BDDC is closely related to FETI-DP \cite{Farhat-2000-SDP}. Though the methods
are quite different, they can be built largely from the same components, and
the eigenvalues of the preconditioned problem (other than the eigenvalue equal
to one) in BDDC and FETI-DP are the same \cite{Mandel-2005-ATP}. See also
\cite{Brenner-2007-BFW,Li-2006-FBB,Mandel-2007-BFM} for simplified proofs.
Thus the performance of BDDC and FETI-DP is essentially identical, and results
for one method apply immediately to the other.

The frontal solver was used to implement a limited variant of BDDC in
\cite{Burda-2007-BMS,Sistek-2007-DEP}. The implementation takes advantage of the
existing integration of the frontal solver into the finite element methodology
and of its implementation of constraints, which is well-suited for BDDC.
However, the frontal solver treats naturally only point constraints, while an
efficient BDDC method in three dimensions requires constraints on averages.
This fact was first observed for FETI-DP experimentally in
\cite{Farhat-2000-SDP}, and theoretically in \cite{Klawonn-2002-DPF}, but
these observations apply to BDDC as well because of the equivalence between
the methods.

In this paper, we extend the previous implementation of BDDC by the frontal
solver to constraints on averages and apply the method to a problem in
biomechanics. We also compare the performance of the mehod with adding
averages and with additional point constraints. The implementation relies on
the separation of point constraints and enforcing the rest by Lagrange
multipliers, as suggested already in \cite{Dohrmann-2003-PSC}. One new aspect
of the present approach is the use of reactions, which come naturally from the
frontal solver, to avoid custom coding.

\section{Mathematical formulation of BDDC}

\label{sec:math}

Consider the problem in a variational form
\begin{equation}
a(u,v)=\langle f,v\rangle\quad\forall\,v\in V\,, \label{eq:var}
\end{equation}
where $V$ is a finite element space of $\mathbb{R}^{3}$-valued piecewise
polynomial functions $v$ continuous on a given domain $\Omega\subset
\mathbb{R}^{3}$, satisfying homogeneous Dirichlet boundary conditions, and
\begin{equation}
a(u,v)=\int_{\Omega}\,(\lambda\,\textup{div}\,u\,\textup{div}\,v+\frac{1}
{2}\,\mu\,(\nabla u+\nabla^{T}u):(\nabla v+\nabla^{T}v)). \label{eq:a}
\end{equation}
Here $\lambda$ and $\mu$ are the first, and the second Lam\' e's parameter,
respectively. Solution $u \in V$ represents the vector field of displacement.
It is known that $a(u,v)$ is a symmetric positive definite bilinear form on
$V$. An equivalent formulation of (\ref{eq:var}) is to find a solution $u$ to
a linear system
\begin{equation}
\label{eq:A-problem}Au=f,
\end{equation}
where $A=(a_{ij})$ is the stiffness matrix computed as $a_{ij}=a(\phi_{i}
,\phi_{j})$, where $\{\phi_{i}\}$ is a finite element basis of $V$,
corresponding to set of unknowns, also called degrees of freedom, defined as
values of displacement at the nodes of a~given triangulation of the domain.
The domain $\Omega$ is decomposed into nonoverlapping subdomains $\Omega_{i}$,
$i=1,\dots N$, also called substructures. Unknowns common to at least two
subdomains are called boundary unknowns and the union of all boundary unknowns
is called the interface $\Gamma$.

The first step is the reduction of the problem to the interface. This is quite
standard and described in the literature, e.g., \cite{Toselli-2005-DDM}. The space $V$
is decomposed as the $a$-orthogonal direct sum $V=V_{1}\oplus\dots\oplus
V_{N}\oplus V_{\Gamma}$, where $V_{i}$ is the space of all functions from $V$
with nonzero values only inside $\Omega_{i}$ (in particular, they are zero on
$\Gamma$), and $V_{\Gamma}$ is the $a$-orthogonal complement of all spaces
$V_{i}$; $V_{\Gamma}=\{v\in V:a(v,w)=0\,\,\forall w\in V_{i},\,i=1,\dots N\}$.
Functions from $V_{\Gamma}$ are fully determined by their values at unknowns
on $\Gamma$ and the discrete harmonic condition that they have minimal energy
on every subdomain (i.e. solve the system with zero right hand side in
corresponding equations). They are represented in the computation by their
values on the interface $\Gamma$. Solution $u$ may be split into the sum of
interior solutions $u_{o}=\sum_{i}^{N}u_{i}$, $u_{i}\in V_{i}$, and
$u_{\Gamma}\in V_{\Gamma}$. Then problem (\ref{eq:A-problem}) may be rewritten
as
\begin{equation}
A(u_{\Gamma}+u_{o})=f. \label{eq:A-problem2}
\end{equation}
Let us now write problem (\ref{eq:A-problem2}) in the block form, with the
first block corresponding to unknowns in subdomain interiors, and the second
block corresponding to unknowns at the interface,
\begin{equation}
\left[
\begin{array}
[c]{cc}
A_{11} & A_{12}\\
A_{21} & A_{22}
\end{array}
\right]  \left[
\begin{array}
[c]{c}
u_{\Gamma1}+u_{o1}\\
u_{\Gamma2}+u_{o2}
\end{array}
\right]  =\left[
\begin{array}
[c]{c}
f_{1}\\
f_{2}
\end{array}
\right]  , \label{eq:A-block}
\end{equation}
with $u_{o2}=0$. Using the fact that functions from $V_{\Gamma}$ are energy
orthogonal to interior functions (so it holds $A_{11} u_{\Gamma1} + A_{12}
u_{\Gamma2} = 0$) and eliminating the variable $u_{o1}$, we obtain that
(\ref{eq:A-block}) is equivalent to
\begin{equation}
A_{11}u_{o1}=f_{1}, \label{eq:A-interior}
\end{equation}
\begin{equation}
\left[
\begin{array}
[c]{cc}
A_{11} & A_{12}\\
A_{21} & A_{22}
\end{array}
\right]  \left[
\begin{array}
[c]{c}
u_{\Gamma1}\\
u_{\Gamma2}
\end{array}
\right]  =\left[
\begin{array}
[c]{c}
0\\
f_{2}-A_{21}u_{o1}
\end{array}
\right]  , \label{eq:A-block-interface2}
\end{equation}
and the solution is obtained as $u=u_{\Gamma}+u_{o}$. Problem
(\ref{eq:A-block-interface2}) is equivalent to the problem
\begin{equation}
\left[
\begin{array}
[c]{cc}
A_{11} & A_{12}\\
0 & S
\end{array}
\right]  \left[
\begin{array}
[c]{c}
u_{\Gamma1}\\
u_{\Gamma2}
\end{array}
\right]  =\left[
\begin{array}
[c]{c}
0\\
g_{2}
\end{array}
\right]  , \label{eq:S-block-interface2}
\end{equation}
where $S$ is the \emph{Schur complement} with respect to interface:
$S=A_{22}-A_{21}A_{11}^{-1}A_{12}$, and $g_{2}$ is the \emph{condensed right
hand side} $g_{2}=f_{2}-A_{21}u_{o1} = f_{2}-A_{21}A_{11}^{-1}f_{1}$. Problem
(\ref{eq:S-block-interface2}) can be split into two problems
\begin{equation}
A_{11}u_{\Gamma1}=-A_{12} u_{\Gamma2}, \label{eq:S-interior}
\end{equation}
\begin{equation}
Su_{\Gamma2}=g_{2}. \label{eq:S-problem}
\end{equation}
Since $A_{11}$ has a block diagonal structure, the solution to
(\ref{eq:A-interior}) may be found in parallel and similarly the solution to
(\ref{eq:S-interior}).

The BDDC method is a particular kind of preconditioner for the reduced problem
(\ref{eq:S-problem}). The main idea of the BDDC preconditioner in an abstract
form \cite{Mandel-2007-BFM} is to construct an auxiliary finite dimensional
space $\widetilde{W}$ such that $V_{\Gamma}\subset\widetilde{W}$ and extend
the bilinear form $a\left(  \cdot,\cdot\right)  $ to a form $\widetilde
{a}\left(  \cdot,\cdot\right)  $ defined on $\widetilde{W}\times\widetilde{W}$
and such that solving the variational problem (\ref{eq:var}) with
$\widetilde{a}\left(  \cdot,\cdot\right)  $ in place of $a\left(  \cdot
,\cdot\right)  $ is cheaper and can be split into independent computations
done in parallel. Then the solution projected to $V_{\Gamma}$ is used for the
preconditioning of $S$. Specifically, let $E:\widetilde{W}\rightarrow
V_{\Gamma}$ be a given projection of $\widetilde{W}$ onto $V_{\Gamma}$, and
$r_{2}={g}_{2}-Su_{\Gamma2}$ the residual in a PCG iteration. Then the output
of the BDDC preconditioner is the part $v_{2}$ of $v=Ew$, where
\begin{equation}
w\in\widetilde{W}:\widetilde{a}\left(  w,z\right)  =\langle r,Ez\rangle
\quad\forall z\in\widetilde{W} , \label{eq:bddc}
\end{equation}
\begin{equation}
r =\left[
\begin{array}
[c]{c}
r_{1}\\
r_{2}
\end{array}
\right]  \, , \quad\quad v =\left[
\begin{array}
[c]{c}
v_{1}\\
v_{2}
\end{array}
\right] \nonumber
\end{equation}
In terms of operators, $v=E\widetilde{S}^{-1}E^{T}r$, where $\widetilde{S}$ is
the operator on $V_{\Gamma}$ associated with the bilinear form $\widetilde{a}$
(but not computed explicitly as a matrix).

Note that while the residual $r_{2}$ in the PCG method applied to the reduced
problem is given at the interface only, the residual in (\ref{eq:bddc}) has
the dimension of all unknowns on the subdomain. This is corrected naturally by
extending the residual $r_{2}$ to subdomain interiors by zeros (setting
$r_{1}=0$), which is required by the condition that the solution $v$ is
discrete harmonic inside subdomain. Similarly, only interface values $v_{2}$
of $v$ are used in further PCG computation. Such approach is equivalent to
computing with explicit Schur complements.

The choice of the space $\widetilde{W}$ and the projection $E$ determines a
particular instance of BDDC \cite{Dohrmann-2003-PSC,Mandel-2007-BFM}. All
functions from $V_{\Gamma}$ are continuous on the domain $\Omega$. In order to
design the space $\widetilde{W}$, we relax the continuity on the interface
$\Gamma$. On $\Gamma$, we select \textit{coarse degrees of freedom} and define
$\widetilde{W}$ as the space of finite element functions with minimal energy
on every subdomain, continuous across $\Gamma$ only at coarse degrees of
freedom. The coarse degrees of freedom can be of two basic types -- explicit
unknowns (called \emph{coarse unknowns}) at selected nodes (called
\emph{corners}), and averages over larger sets of nodes (subdomain
\emph{faces} or \emph{edges}). The continuity condition then means that the
values at the corresponding corners, resp. averages, on neighbouring
subdomains coincide. The bilinear form $a\left(  \cdot,\cdot\right)  $ from
(\ref{eq:a}) is extended to $\widetilde{a}\left(  \cdot,\cdot\right)  $ on
$\widetilde{W}\times\widetilde{W}$ by integrating (\ref{eq:a}) over the
subdomains $\Omega_{i}$ separately and adding the results.

The projection $E:\widetilde{W}\rightarrow V_{\Gamma}$ is defined at unknowns
on the interface $\Gamma$ (the $u_{\Gamma2}$ part) as a weighted average of
values from different subdomains and thus resulting in function continuous
across the interface. These averaged values on $\Gamma$ determine the
projection $E$, because values inside subdomains (the $u_{\Gamma1}$ part) are
then obtained by the solutions of local subdomain problems
(\ref{eq:S-interior}) to make the averaged function discrete harmonic. To
assure good performance regardless of different stiffness of the subdomains
\cite{Mandel-2003-CBD}, the weights are chosen proportional to the
corresponding diagonal entries of the subdomain stiffness matrices. The
transposed projection $E^{T}$ is used for distribution of the residual $r_{2}$
among neighbouring subdomains and represents the decomposition of unity at
unknowns on interface.

The decomposition into subspaces used to derive the problem with Schur
complement (\ref{eq:S-problem}) is now repeated for space $\widetilde{W}$,
with the coarse degrees of freedom playing the role of interface unknowns and
$\widetilde{a}(\cdot,\cdot)$ the role of $a(\cdot,\cdot)$. Namely, space
$\widetilde{W}$ is decomposed as $\widetilde{a}$-orthogonal direct sum
$\widetilde{W}=\widetilde{W}_{1}\oplus\dots\oplus\widetilde{W}_{N}
\oplus\widetilde{W}_{C}$, where $\widetilde{W}_{i}$ is the space of functions
with nonzero values only in $\Omega_{i}$ outside coarse degrees of freedom
(they have zero values at corners, they are generally not continuous at other
unknowns on $\Gamma$, and they have zero averages) and $\widetilde{W}_{C}$ is
the \textit{coarse space}, defined as the $\widetilde{a}$-orthogonal
complement of all spaces $\widetilde{W}_{i}$: $\widetilde{W}_{C}
=\{v\in\widetilde{W}:\widetilde{a}(v,w)=0\ \forall w\in\widetilde{W}
_{i},\ i=1,\dots N\}$. Functions from $\widetilde{W}_{C}$ are fully determined
by their values at coarse degrees of freedom (where they are continuous) and
have minimal energy. Thus, they are generally discontinuous across $\Gamma$
outside the corners. The solution $w\in\widetilde{W}$ from (\ref{eq:bddc}) is
now split accordingly as $w=w_{C}+\sum_{i=1}^{N}w_{i},$ where $w_{C}$,
determined by
\begin{equation}
w_{C}\in\widetilde{W}_{C}:\widetilde{a}\left(  w_{C},v\right)  =\langle
r,Ev\rangle\quad\forall v\in\widetilde{W}_{C}, \label{eq:coarse-problem}
\end{equation}
is called the \emph{coarse correction}, and $w_{i}$, determined by
\begin{equation}
w_{i}\in\widetilde{W}_{i}:\widetilde{a}\left(  w_{i},v\right)  =\langle
r,Ev\rangle\quad\forall v\in\widetilde{W}_{i}, \label{eq:subdomain-problem}
\end{equation}
is the \emph{substructure correction} from $\Omega_{i}$, $i=1,\dots N$.

Let us now rewrite the BDDC preconditioner in terms of matrices, following
\cite{Dohrmann-2003-PSC}. Problem (\ref{eq:subdomain-problem}) is formulated
in a saddle point form as
\begin{equation}
\left[
\begin{array}
[c]{cc}
K_{i} & \overline{C}_{i}^{T}\\
\overline{C}_{i} & 0
\end{array}
\right]  \left[
\begin{array}
[c]{c}
w_{i}\\
\mu_{i}
\end{array}
\right]  =\left[
\begin{array}
[c]{c}
r_{i}\\
0
\end{array}
\right]  , \label{eq:subdomain_matrix1}
\end{equation}
where $K_{i}$ denotes the substructure local stiffness matrix, obtained by the
subassembly of element matrices only of elements in substructure $i$, matrix
$\overline{C}_{i}$ represents constraints on subdomain, that enforce zero
values of coarse degrees of freedom, $\mu_{i}$ is vector of Lagrange
multipliers, and $r_{i}$ is the weighted residual $E^{T} r$ restricted to
subdomain $i$.

Matrix $K_{i}$ is singular for floating subdomains (subdomains not touching
Dirichlet boundary conditions), while the augmented matrix of problem
(\ref{eq:subdomain_matrix1}) is regular and may be factorized. Matrix
$\overline{C}_{i}$ contains both constraints enforcing continuity across
corners (single point continuity), and constraints enforcing equality of
averages over edges and faces of subdomains. The former type corresponds to
just one nonzero entry equal to 1 on a row of $\overline{C}_{i}$, while the
latter leads to several nonzero entries on a row. This structure will be
exploited in the following section.

Problem (\ref{eq:subdomain_matrix1}) is solved in each iteration of the PCG
method to find the correction from substructure $i$. However, the matrix of
(\ref{eq:subdomain_matrix1}) is used prior the whole iteration process to
construct the \emph{local subdomain matrix of the coarse problem}. First, the
\emph{coarse basis functions} are found independently for each subdomain as
the solution to
\begin{equation}
\left[
\begin{array}
[c]{cc}
K_{i} & \overline{C}_{i}^{T}\\
\overline{C}_{i} & 0
\end{array}
\right]  \left[
\begin{array}
[c]{c}
\psi_{i}\\
\lambda_{i}
\end{array}
\right]  =\left[
\begin{array}
[c]{c}
0\\
I
\end{array}
\right]  . \label{eq:subdomain_matrix2}
\end{equation}
This is a problem with multiple right hand sides, where $\psi_{i}$ is a matrix
of coarse basis functions with several columns, each corresponding to one
coarse degree of freedom on subdomain. These functions are given by values
equal to 0 at all coarse degrees of freedom except one, where they have value
equal to 1, and they have minimal energy on subdomain outside coarse degrees
of freedom. The identity block $I$ has the dimension of the number of
constraints on the subdomain.

Once $\psi_{i}$ is known, the \emph{subdomain coarse matrix} $K_{Ci}$ is
constructed as
\begin{equation}
\label{eq:coarse_matrix1}K_{Ci} = \psi_{i}^{T} K_{i} \psi_{i} .
\end{equation}

Matrices $K_{Ci}$ are then assembled to form the global \emph{coarse matrix}
$A_{C}$. This procedure is same as the standard process of assembly in finite
element solution, with subdomains playing the role of elements, coarse degrees
of freedom on subdomain representing degrees of freedom on element, and matrix
$K_{Ci}$ representing the element stiffness matrix.

Problem (\ref{eq:coarse-problem}) is now
\begin{equation}
A_{C}\mathbf{w}_{C}=r_{C}, \label{eq:coarse_matrix2}
\end{equation}
where $r_{C}$ is the \emph{global coarse residual} obtained by the assembly of
the subdomain contributions of the form $r_{Ci}=\psi_{i}^{T}r_{i}$.

The coarse solution $\mathbf{w}_{C}$ has the dimension of the number of all
coarse degrees of freedom. So, to add the correction to subdomain problems, we
first have to restrict it to coarse degrees of freedom on each subdomain and
to interpolate it to the whole subdomain by $w_{Ci}=\psi_{i}\mathbf{w}_{Ci}$.
By extending $w_{Ci}$ and $w_{i}$ by zero to other subdomains, these can be
summed over the subdomains to form the final vector $w$. Finally, the
preconditioned residual is obtained as $v=Ew$.

It is worth noticing that in the case of no constraints on averages, i.e.
using only coarse unknowns for the definition of the coarse space, matrix
$A_{C}$ of problem (\ref{eq:coarse_matrix2}) is simply the Schur complement of
matrix $A$ with respect to coarse unknowns. This fact was pointed out in
\cite{Li-2006-FBB}. If additional degrees of freedom are added for averages,
they correspond to new explicit unknowns in $\mathbf{w}_{Ci}$.

Obviously, several mapping operators among various spaces are needed in the
implementation, defining embedding of subdomains into global problem, local
subdomain coarse problem into global coarse problem etc. We have circumvented
their mathematical definition by words for the sake of brevity, while we refer
to \cite{Dohrmann-2003-PSC,Mandel-2003-CBD} for rigorous definitions of these operators.

\section{BDDC implementation based on frontal solver}

The frontal solver implements the solution of a square linear system with some
of the variables having prescribed values. Equations that correspond to the
fixed variables are omitted and the values of these variables are substituted
into the solution vector directly. The output of the solver consists of the
solution and the resulting imbalance in the equations, called reaction forces.
More precisely, consider a block decomposition of the vector of unknowns $x$
with the second block consisting of all fixed variables, and write a system
matrix $A$ with the same block decomposition (here, the decomposition is
different from the one in Section \ref{sec:math}). Then on exit from the
frontal solver,
\begin{equation}
\left[
\begin{array}
[c]{cc}
A_{11} & A_{12}\\
A_{21} & A_{22}
\end{array}
\right]  \left[
\begin{array}
[c]{c}
x_{1}\\
x_{2}
\end{array}
\right]  =\left[
\begin{array}
[c]{c}
f_{1}\\
f_{2}
\end{array}
\right]  +\left[
\begin{array}
[c]{c}
0\\
r_{2}
\end{array}
\right]  , \label{eq:frontal}
\end{equation}
where fixed variable values $x_{2}$ and the load vectors $f_{1}$ and $f_{2}$
are the inputs, while the solution $x_{1}$ and the reaction $r_{2}$ are the outputs.
Stiffness matrices of elements are input instead of the whole matrix,
and their assembly is done simultaneously with the factorization inside the frontal solver.

The key idea of this section is to split the constraints in matrix
$\overline{C}_{i}$ and to handle them in different ways. Those enforcing zero
values at corners will be enforced as fixed variables, while the remaining
constraints, corresponding to averages and denoted $C_{i}$, will be still
enforced using Lagrange multipliers.the

In the rest of this section, we drop the subdomain subscript $_{i}$ and we
write subdomain vectors $w$ in the block form with the second block consisting
of unknowns that are involved in coarse degrees of freedom (i.e. coarse
unknowns), denoted by the subscript $_{c}$, and the first block consisting of
the remaining degrees of freedom, denoted by the subscript $_{f}$. The vector
of the coarse degrees of freedom given by averages is written as $Cw$, where
each row of $C$ contains the coefficients of the average that makes that
degree of freedom. Then subdomain vectors $w\in\widetilde{W}$ are
characterized by $w_{c}=0$, $Cw=0$. Assume that $C=\left[  C_{f}
\ \ C_{c}\right]  $, with $C_{c}=0$, that is, the averages do not involve
single variable coarse degrees of freedom; then $Cw=C_{f}w_{f}$. The subdomain
stiffness matrix $K$ is singular for floating subdomains, but the block
$K_{ff}$ is nonsingular if there are enough corners to eliminate the rigid
body motions, which will be assumed.

We now show how to solve (\ref{eq:subdomain_matrix1}) --
(\ref{eq:coarse_matrix2}) using the frontal solver. In the case when there are
no averages as coarse degrees of freedom, we recover the previous method from
\cite{Burda-2007-BMS,Sistek-2007-DEP}.

The local substructure problems (\ref{eq:subdomain-problem}) are written in
the frontal solver form (\ref{eq:frontal}) as
\begin{equation}
\left[
\begin{array}
[c]{ccc}
K_{ff} & K_{fc} & C_{f}^{T}\\
K_{cf} & K_{cc} & 0\\
C_{f} & 0 & 0
\end{array}
\right]  \left[
\begin{array}
[c]{c}
w_{f}\\
w_{c}\\
\mu
\end{array}
\right]  =\left[
\begin{array}
[c]{c}
r\\
0\\
0
\end{array}
\right]  +\left[
\begin{array}
[c]{c}
0\\
Rea\\
0
\end{array}
\right]  , \label{eq:subdomain-frontal}
\end{equation}
where $w_{c}=0$, $r$ is the part in the $_{f}$ block of the residual in the
PCG method distributed to the substructures by the operator $E^{T}$, and $Rea$
is the reaction. The constraint $w_{c}=0$ is enforced by marking the $w_{c}$
unknowns as fixed, while the remaining constraints $C_{f}w_{f}=0$ are enforced
via the Lagrange multiplier $\mu$. Using the fact that $w_{c}=0$, we get from
(\ref{eq:subdomain-frontal}) that
\begin{align}
K_{ff}w_{f}  &  =-C_{f}^{T}\mu+r,\label{eq:subs-sys-1}\\
K_{cf}w_{f}  &  =Rea,\label{eq:subs-sys-2}\\
C_{f}w_{f}  &  =0. \label{eq:subs-sys-3}
\end{align}
From (\ref{eq:subs-sys-1}), $w_{f}=K_{ff}^{-1}\left(  -C_{f}^{T}\mu+r\right)
$. Now substituting $w_{f}$ into (\ref{eq:subs-sys-3}), we get the problem for
Lagrange multiplier $\mu$,
\begin{equation}
C_{f}K_{ff}^{-1}C_{f}^{T}\mu=C_{f}K_{ff}^{-1}r. \label{eq:dual_sub}
\end{equation}

The matrix $C_{f}K_{ff}^{-1}C_{f}^{T}$ is dense but small, with the order
equal to the number of averages on the subdomain, and it is constructed by
solving the system $K_{ff}U=C_{f}^{T}$ with multiple right hand sides by the
frontal solver and then the multiplication $C_{f}U$. After solving problem
(\ref{eq:dual_sub}), we substitute for $\mu$ in (\ref{eq:subs-sys-1}) and find
$w_{f}$ from
\begin{equation}
\left[
\begin{array}
[c]{cc}
K_{ff} & K_{fc}\\
K_{cf} & K_{cc}
\end{array}
\right]  \left[
\begin{array}
[c]{c}
w_{f}\\
w_{c}
\end{array}
\right]  =\left[
\begin{array}
[c]{c}
r - C_{f}^{T} \mu\\
0
\end{array}
\right]  +\left[
\begin{array}
[c]{c}
0\\
Rea
\end{array}
\right]  \label{eq:subdomain-frontal2}
\end{equation}
by the frontal solver, considering $w_{c}=0$ fixed. The factorization in the
frontal solver for (\ref{eq:subdomain-frontal2}) and the factorization of the
matrix $C_{f}K_{ff}^{-1}C_{f}^{T}$ for (\ref{eq:dual_sub}) need to be computed
only once in the setup phase.

The coarse problem (\ref{eq:coarse_matrix2}) is solved by the frontal solver
just like a~finite element problem, with the subdomains playing the role of
elements. It only remains to specify the basis functions of $\widetilde{W}
_{C}$ on the subdomain from (\ref{eq:subdomain_matrix2}), and compute the
local subdomain coarse matrix (\ref{eq:coarse_matrix1}) efficiently. Denote by
$\psi^{c}$ the matrix whose colums are coarse basis functions associated with
the coarse unknowns at corners, and $\psi^{avg}$ the matrix made out of the
coarse basis functions associated with averages. To find the coarse basis
functions, we proceed similarly as in (\ref{eq:subdomain-frontal}) and write
the equations for the coarse basis functions in the frontal solver form, now
with multiple right-hand sides,
\begin{equation}
\left[
\begin{array}
[c]{ccc}
K_{ff} & K_{fc} & C_{f}^{T}\\
K_{cf} & K_{cc} & 0\\
C_{f} & 0 & 0
\end{array}
\right]  \left[
\begin{array}
[c]{cc}
\psi_{f}^{c} & \psi_{f}^{avg}\\
I & 0\\
\lambda^{c} & \lambda^{avg}
\end{array}
\right]  =\left[
\begin{array}
[c]{cc}
0 & 0\\
0 & 0\\
0 & I
\end{array}
\right]  +\left[
\begin{array}
[c]{cc}
0 & 0\\
Rea^{c} & Rea^{avg}\\
0 & 0
\end{array}
\right]  , \label{eq:block1}
\end{equation}
where $Rea^{c}$ and $Rea^{avg}$ are matrices of reactions. Denote $\psi
_{f}=\left[
\begin{array}
[c]{cc}
\psi_{f}^{c} & \psi_{f}^{avg}
\end{array}
\right]  $, $\psi_{c}=\left[
\begin{array}
[c]{cc}
\psi_{c}^{c} & \psi_{c}^{avg}
\end{array}
\right]  =\left[
\begin{array}
[c]{cc}
I & 0
\end{array}
\right]  $, $\lambda=\left[
\begin{array}
[c]{cc}
\lambda^{c} & \lambda^{avg}
\end{array}
\right]  $, $Rea=\left[
\begin{array}
[c]{cc}
Rea^{c} & Rea^{avg}
\end{array}
\right]  $, and $R=\left[
\begin{array}
[c]{cc}
0 & I
\end{array}
\right]  $. Then (\ref{eq:block1}) becomes
\begin{align}
K_{ff}\psi_{f}+K_{fc}\psi_{c}  &  =-C_{f}^{T}\lambda,\label{eq:coarse-sys-1}\\
K_{cf}\psi_{f}+K_{cc}\psi_{c}  &  =Rea,\label{eq:coarse-sys-2}\\
C_{f}\psi_{f}  &  =R. \label{eq:coarse-sys-3}
\end{align}
From (\ref{eq:coarse-sys-1}), we get $\psi_{f}=-K_{ff}^{-1}\left(  K_{fc}
\psi_{c}+C_{f}^{T}\lambda\right)  $. Substituting $\psi_{f}$ into
(\ref{eq:coarse-sys-3}), we derive the problem for Lagrange multipliers
\begin{equation}
C_{f}K_{ff}^{-1}C_{f}^{T}\lambda=-\left(  R+C_{f}K_{ff}^{-1}K_{fc}\psi
_{c}\right)  , \label{eq:dual-sys-1}
\end{equation}
which is solved for $\lambda$ by solving the system (\ref{eq:dual-sys-1}) for
multiple right hand sides. Since $\psi_{c}$ is known, we can use the frontal
solver to solve (\ref{eq:coarse-sys-1})-(\ref{eq:coarse-sys-2}) to find
$\psi_{f}$ and $Rea$:

\begin{equation}
\left[
\begin{array}
[c]{cc}
K_{ff} & K_{fc}\\
K_{cf} & K_{cc}
\end{array}
\right]  \left[
\begin{array}
[c]{c}
\psi_{f}\\
\psi_{c}
\end{array}
\right]  =\left[
\begin{array}
[c]{c}
- C_{f}^{T} \lambda\\
0
\end{array}
\right]  +\left[
\begin{array}
[c]{c}
0\\
Rea
\end{array}
\right]  , \label{eq:subdomain-frontal3}
\end{equation}

considering $\psi_{c}$ fixed. Finally, we construct the local coarse matrix
corresponding to the subdomain as
\begin{equation}
K_{C}=\psi^{T}K\psi=\psi^{T}\left[
\begin{array}
[c]{c}
-C_{f}^{T}\lambda\\
Rea
\end{array}
\right]  =\left[
\begin{array}
[c]{cc}
\psi_{f}^{T} & \psi_{c}^{T}
\end{array}
\right]  \left[
\begin{array}
[c]{c}
-C_{f}^{T}\lambda\\
Rea
\end{array}
\right]  =-\psi_{f}^{T}C_{f}^{T}\lambda+\left[
\begin{array}
[c]{c}
I\\
0
\end{array}
\right]  Rea, \label{eq:coarse-create}
\end{equation}
where $\psi=\left[
\begin{array}
[c]{cc}
\psi^{c} & \psi^{avg}
\end{array}
\right]  $.

At the end of the setup phase, the matrix of coarse problem is factored by the
frontal solver, using subdomain coarse matrices as input.

\section{The algorithm}

The setup starts with two types of factorization of the system matrix by
frontal algorithm that differ only in the set of fixed variables. In the first
factorization, all interface unknowns are prescribed as fixed, for the
solution of problems (\ref{eq:A-block}) and (\ref{eq:A-block-interface2}). In
the second factorization, only the coarse unknowns are fixed, for the solution
of problems (\ref{eq:subdomain-frontal2}) and (\ref{eq:subdomain-frontal3}).
Both factorizations are done subdomain by subdomain, so the interface unknowns
are represented by their own instance in different subdomains. This naturally
leads to parallelization according to subdomains. Then the main steps of the
algorithm are as follows.

(A) Solve (\ref{eq:A-interior}) in parallel on every subdomain for $g_{2}$
as reaction using the frontal solver with $u_{o2}=0$ fixed,
\begin{equation}
\left[
\begin{array}
[c]{cc}
A_{11} & A_{12}\\
A_{21} & A_{22}
\end{array}
\right]  \left[
\begin{array}
[c]{c}
u_{o1}\\
0
\end{array}
\right]  =\left[
\begin{array}
[c]{c}
f_{1}\\
f_{2}
\end{array}
\right]  +\left[
\begin{array}
[c]{c}
0\\
-g_{2}
\end{array}
\right]  \label{eq:A-interior3}
\end{equation}
(B) Solve (\ref{eq:S-problem}) for $u_{\Gamma2}$ using PCG with BDDC
preconditioner (for more details see bellow).\newline(C) Compute the solution
$u$ from (\ref{eq:A-block}) using the frontal solver with $u_{\Gamma2}$ fixed
and $u_{o2}=0$.

Step (B) in detail follows. Before starting cycle of PCG iterations,
compute in advance:

\begin{itemize}
\item Matrix $C_{f}K_{ff}^{-1}C_{f}^{T}$ of (\ref{eq:dual_sub}), resp.
(\ref{eq:dual-sys-1}) and its factorization.

\item Coarse basis functions $\psi_{i}$ of (\ref{eq:subdomain_matrix2}) from
(\ref{eq:subdomain-frontal3}). Multiplier $\lambda$ on the right hand side of
(\ref{eq:subdomain-frontal3}) is computed from (\ref{eq:dual-sys-1}).

\item Local matrices $K_{Ci}$ of (\ref{eq:coarse_matrix1}) from
(\ref{eq:coarse-create}).

\item Factorization of matrix $A_{C}$ of (\ref{eq:coarse_matrix2}) using the
frontal solver with local coarse matrices $K_{Ci}$ as `element' matrices.

\item The first residual $r_{2}$ of (\ref{eq:S-problem}) from the first
approximation of $u_{\Gamma2}$ using the frontal solver on
(\ref{eq:A-block-interface2}) with $u_{\Gamma2}$ fixed,
\begin{equation}
\left[
\begin{array}
[c]{cc}
A_{11} & A_{12}\\
A_{21} & A_{22}
\end{array}
\right]  \left[
\begin{array}
[c]{c}
u_{\Gamma1}\\
u_{\Gamma2}
\end{array}
\right]  =\left[
\begin{array}
[c]{c}
0\\
g_{2}
\end{array}
\right]  +\left[
\begin{array}
[c]{c}
0\\
-r_{2}
\end{array}
\right]  . \label{eq:A-interior4}
\end{equation}
For the popular choice of initial solution $u_{\Gamma 2} = 0$,
this step reduces to setting $r_2 = g_2$
(the discrete harmonic extension of zero function on interface is zero also in the interior of subdomain).

\end{itemize}

Then in every PCG iteration, compute $v_{2}$ from $r_{2}$ in three steps:

\begin{enumerate}
\item Compute $E^{T} r$ by distributing the residual $r_{2}$ on interface
$\Gamma$ among neigbouring subdomains. \newline For every subdomain, compute
$r_{i}$ as restriction of $E^{T} r$ to that subdomain.

\item Compute $w = \widetilde{S}^{-1}E^{T}r$ as sum of all substructure
corrections $w_{i}$ and coarse correction $w_{C}$. These can be computed
in parallel for every $r_{i}$:

\begin{itemize}
\item Substructure correction $w_{i}$ given by (\ref{eq:subdomain_matrix1}) is
computed from (\ref{eq:subdomain-frontal2}) using frontal with $w_{c} = 0$
fixed. Multiplier $\mu$ on the right hand side of (\ref{eq:subdomain-frontal2}
) is computed from (\ref{eq:dual_sub}) and $r$ in both
(\ref{eq:subdomain-frontal2}) and (\ref{eq:dual_sub}) is the part in the
$_{f}$ block of the $r_{i}$.

\item Coarse correction $w_{C}$ is solution of (\ref{eq:coarse_matrix2}).
\end{itemize}

\noindent We are interested in values of $w$ only on interface.

\item Compute the $v_{2}$ part of $v = E w$ at every interface node as
weighted average of values of $w$ at that node (neigbouring subdomains have
generally different values of $w$ at corresponding interface nodes).
\end{enumerate}

Note that in every iteration of PCG, the product $Sp_{2}$ is needed, where
$p_{2}$ is a search direction. This product is computed using the frontal
solver with $p_{2}$ fixed,
\begin{equation}
\left[
\begin{array}
[c]{cc}
A_{11} & A_{12}\\
A_{21} & A_{22}
\end{array}
\right]  \left[
\begin{array}
[c]{c}
p_{1}\\
p_{2}
\end{array}
\right]  =\left[
\begin{array}
[c]{c}
0\\
0
\end{array}
\right]  +\left[
\begin{array}
[c]{c}
0\\
Sp_{2}
\end{array}
\right]  . \label{eq:A-interior5}
\end{equation}
The interior part $p_{1}$ is computed only as a by-product, and it is not used
in the PCG iterations.

\section{Numerical results}

\begin{figure}[p]
\centering  \hbox{ \vbox{ \hsize=70mm \centering
\includegraphics[width=50mm]{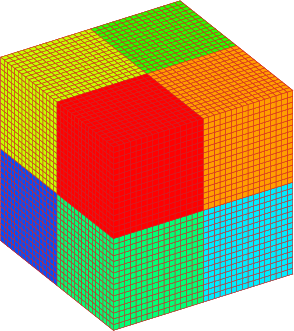}
} \vbox{ \hsize=70mm \centering
\includegraphics[width=50mm]{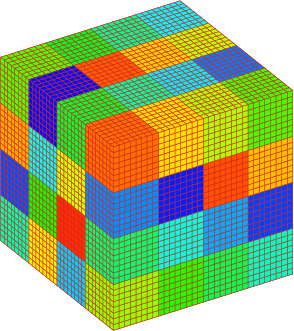}
}
}\caption{Cube problem, division into 8 (left) and 64 (right) subdomains}
\label{fig:cube}
\end{figure}

The implementation was first tested on the problem of unit cube, a classical
test problem of domain decomposition methods. In our case, the cube is made of
steel with Young's modulus $2.1\cdot10^{11}$ Pa and Poisson's ratio $0.3$. The
cube is fixed at one face and loaded by the force of $1,000$ N, acting on one
edge opposite to the fixed face in direction parallel to it and pointing
outwards of the cube. The mesh consists of $32^{3}=32,768$ trilinear elements.
It was uniformly divided into $8$ and $64$ subdomains, resulting in $H/h=16$
and $H/h=8$, respectively (here $H$ denotes the characteristic size of
subdomains and $h$ the characteristic size of elements). These divisions are
presented in Figure \ref{fig:cube}. The interface is initially divided into
$7$ corners, $6$ edges, and $12$ faces in the case of $8$ subdomains, and into
$81$ corners, $108$ edges, and $144$ faces in the case of $64$ subdomains.

All experiments with this problem were computed on 8 1.5~GHz Intel Itanium~2
processors of SGI Altix~4700 computer in CTU Supercomputing Centre, Prague.
The stopping criterion of PCG was chosen as $\Vert r\Vert_{2}/\Vert g\Vert
_{2}<10^{-6}$. In the presented results, an external parallel multifrontal
solver MUMPS \cite{Amestoy-2000-MPD} was used for the factorization and
solution of the coarse problem (\ref{eq:coarse_matrix2}), instead of the
serial frontal solver.

\begin{figure}[p]
\begin{center}
\includegraphics[width=80mm,angle=-90]{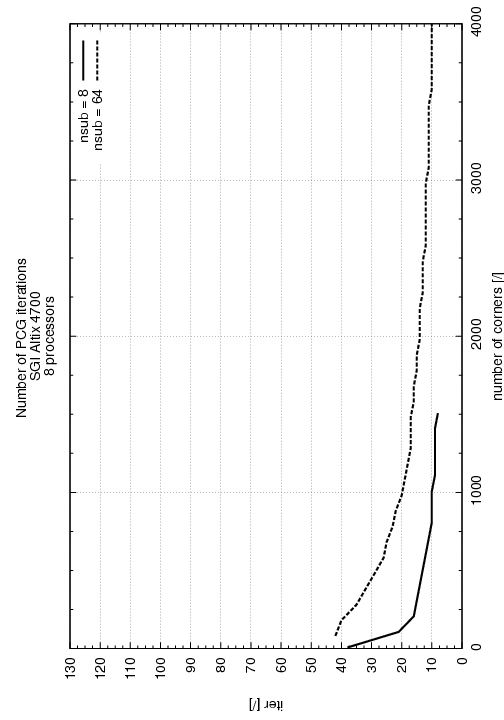}
\end{center}
\caption{Cube problem, number of iterations for adding corners}
\label{fig:cube_iter}
\end{figure}

\begin{figure}[p]
\begin{center}
\includegraphics[width=80mm,angle=-90]{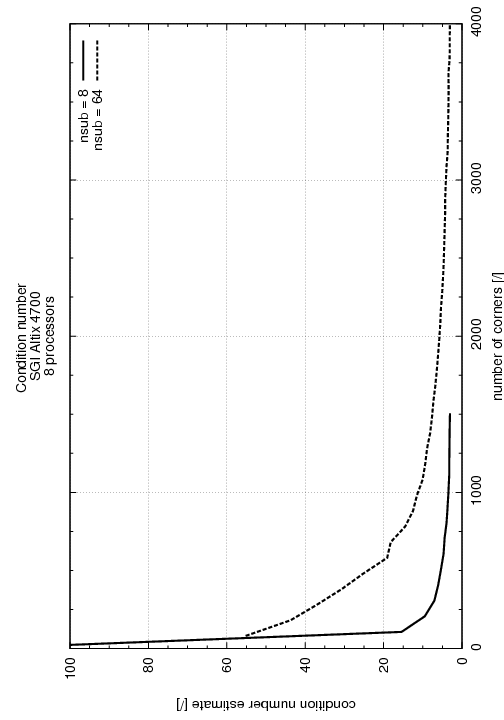}
\end{center}
\caption{Cube problem, condition number for adding corners}
\label{fig:cube_condition_number}
\end{figure}

\begin{figure}[p]
\begin{center}
\includegraphics[width=80mm,angle=-90]{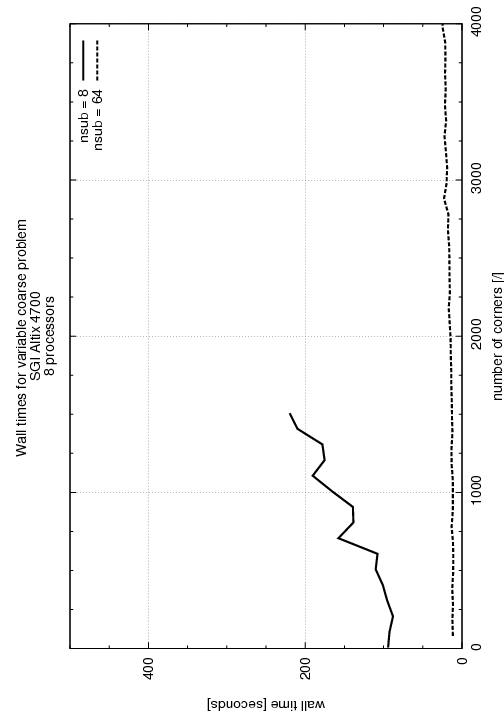}
\end{center}
\caption{Cube problem, wall clock time for adding corners}
\label{fig:cube_wall_time}
\end{figure}

The first experiment compares two ways of enriching the coarse space
$\widetilde{W}_{C}$, namely by adding \emph{point constraints} on randomly
selected variables on the substructure interfaces, i.e. adding more
\textquotedblleft corners\textquotedblright\ (Fig. \ref{fig:cube_iter} --
\ref{fig:cube_wall_time}), and by adding averages to the initial set of
corners (Table \ref{tab:cube8} and Table \ref{tab:cube64}). In the first
column of these tables, no additional averages are considered and only corners
were used in the construction of $\widetilde{W}_{C}$. Then we enforce the
equality of arithmetic averages over all edges, over all faces, and over all
edges and faces, respectively.

Although adding corners leads to an improvement of preconditioner in terms of
the condition number (Fig. \ref{fig:cube_condition_number}) and the number of
iterations (Fig. \ref{fig:cube_iter}), after a slight decrease early on the
total computational time increases (Fig. \ref{fig:cube_wall_time}) due to the
added cost of the setup and factorization of the coarse problem. This effect
is particularly pronounced with $8$ subdomains, where the cost of creating the
coarse matrix dominates, as the frontal solver internally involves
multiplication of large dense matrices to compute reactions:
(\ref{eq:subdomain-frontal3}) takes $O(n_{i}n_{ci}^{2})$ for multiple right
hand sides, where $n_{i}$ is the number of variables and $n_{ci}$ is the
number of coarse variables in subdomain $i$. The problem divided into $64$
subdomains requires much less time than the problem with $8$ subdomains also
due to the fact that the factorization time for subdomain problems grows fast
with subdomain size. Note that for $64$ subdomains, the number of processors
remains the same and each of the $8$ processors handles $8$ subdomains.

\begin{figure}[p]
\begin{center}
\includegraphics[width=100mm]{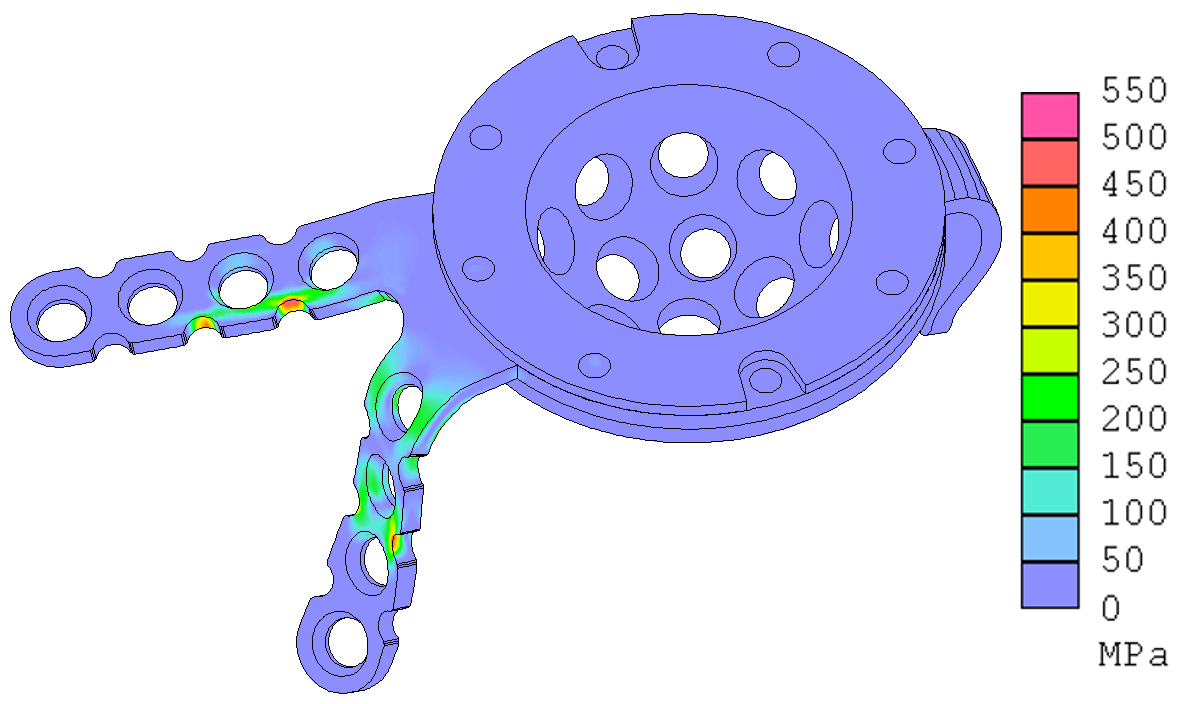}\newline
\end{center}
\caption{Hip joint replacement, von Mises stresses in improved design.}
\label{fig:stress}
\end{figure}

The structural analysis of the replacement of the hip joint construction
loaded by pressure from body weight is an important problem in bioengineering.
The hip replacement consists of several parts made of titanium; here we
consider the central part of the replacement joint. The problem was simplified
to stationary linearized elasticity. The highest stress was reached in the
notches of the holders. In the original design, holders of the hip replacement
had thickness of $2$ mm, which led to maximal von Mises stress about $1,500$
MPa. As the yield point of titanium is about $800$ MPa, the geometry of the
construction had to be modified. The thickness of the holders was increased to
3 mm, radiuses of the notches were increased, and the notches were made
smaller, as in Fig.~\ref{fig:stress}. The maximal von Mises stress on this new
construction was only about 540 MPa, which satisfed the demands for the
strength of the construction \cite{Certikova-2005-SCH,Tuzar-2005-THE}. The mesh
consists of 33,186 quadratic elements resulting in 544,734 unknowns.

The computation needs 400 minutes when using a serial frontal solver on Compaq
Alpha server ES47 at the Institute of Thermomechanics, Academy of Sciences of
the Czech Republic. With 32 subdomains and corner coarse degrees of freedom
only, BDDC on a single Alpha processor took 10 times less, only 40 minutes.

In this paper, we present results for decompositions into 16 and 32 subdomains
obtained by the METIS graph partitioner \cite{Karypis-1998-FHQ}. In the case
of 16 subdomains, the interface topology leads to 35 corners, 12 edges, and 35
faces, and in the case of 32 subdomains, to 57 corners, 12 edges, and 66
faces. All presented results were obtained on 16 processors of SGI Altix~4700.
Again, the stopping criterion of PCG was chosen as $\|r\|_{2}/\|g\|_{2} <
10^{-6}$. Also these results were obtained using MUMPS solver for the coarse
problem solution.

\begin{figure}[p]
\begin{center}
\includegraphics[width=80mm,angle=-90]{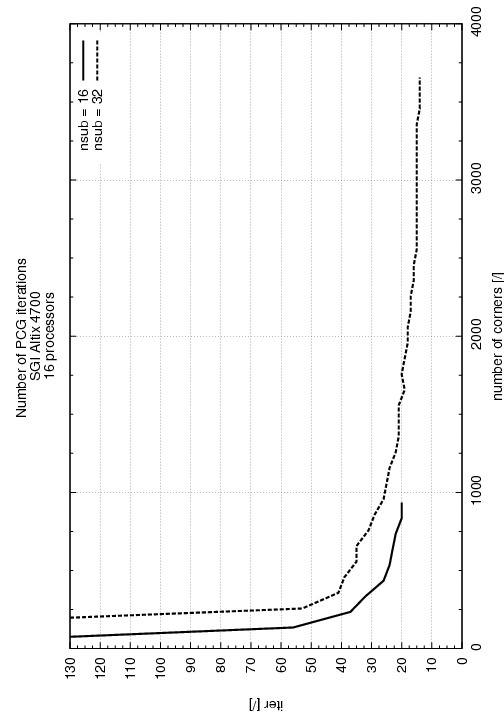}
\end{center}
\caption{Hip joint replacement, number of iterations for adding corners}
\label{fig:hip_iter}
\end{figure}

\begin{figure}[p]
\begin{center}
\includegraphics[width=80mm,angle=-90]{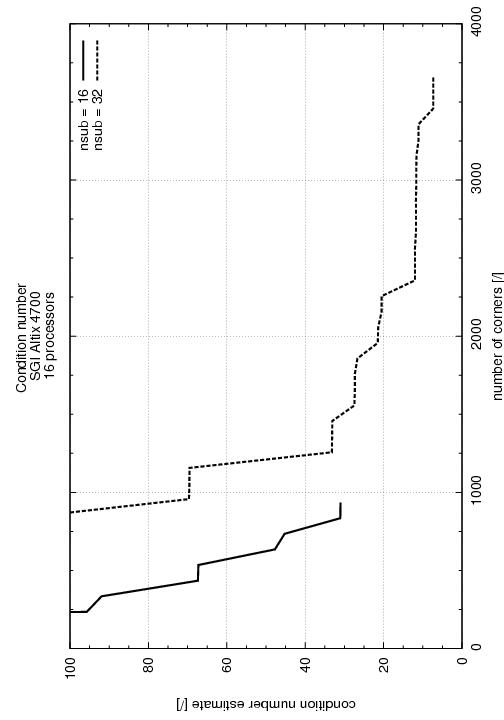}
\end{center}
\caption{Hip joint replacement, condition number for adding corners}
\label{fig:hip_condition_number}
\end{figure}

\begin{figure}[p]
\begin{center}
\includegraphics[width=80mm,angle=-90]{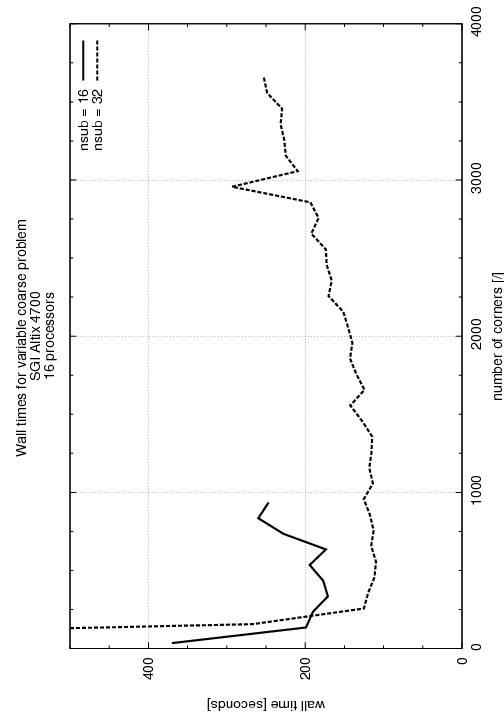}
\end{center}
\caption{Hip joint replacement, wall clock time for adding corners}
\label{fig:hip_wall_time}
\end{figure}

We again investigate the two approaches to coarse space enrichment. The
results of random addition of corners are presented in Figures
\ref{fig:hip_iter} -- \ref{fig:hip_wall_time}. Unlike in the case of the cube,
a substantial decrease of the total time is achieved. The optimal number of
additional corners depends on the division into subdomains.

Results of adding averages are summarized in Table \ref{tab:hip16basic} for
the case of 16 subdomains, and in Table \ref{tab:hip32basic} for the case of
32 subdomains, both with the initial set of corners.

The last experiment combines both approaches: we add averages to the optimal
size of the set of corners, determined from Figure \ref{fig:hip_wall_time} as
${\sim}335$ for the problem with $16$ subdomains, and ${\sim}557$ for the
problem with $32$ subdomains. These results are presented in Tables
\ref{tab:hip16optimal} and \ref{tab:hip32optimal}. We can see that this
synergy can lead to the lowest overall time of the computation.

We observe that while the implementation of averages leads to a negligible
increase in the computational cost of the factorizations, it considerably
improves the condition number, and thus reduces the overall time of solution.
Also, the decomposition into 32 subdomains leads to significantly lower
computational times than the division into 16 subdomains. Note that although
the initial set of corners leads to non-singular local matrices and the coarse
matrix and so successful setup of the preconditioner, the iterations do not
converge in this case.

\section{Conclusion}

We have presented an application of a standard frontal solver within the
iterative substructuring method BDDC. The method was applied to a
biomechanical stress analysis problem. The numerical results show that the
improvement of preconditioning by additional constraints is significant and
can lead to a considerable savings of computational time, while the additional
cost is negligible.

For a model problem (cube), constraints by averages on edges and faces are
required for good performance, as predicted by the theory
\cite{Klawonn-2002-DPF,Mandel-2005-ATP}, and additional point constraints
(i.e., corners) are not productive. However, for the hip replacement problem
(which is far from a regularly decomposed cube), additional point constraints
result in significantly lower total computational time than the averages, and
the best result is obtained by combining both the added point constraints and
the averages.

For large problems and a large number of processors, load balancing will be
essential \cite{Medek-2007-SLB}.

\section*{Acknowledgement}

This research has been supported by the Czech Science Foundation under grant
106/08/0403 and by the U.S. National Science Foundation under grant
DMS-0713876. It has also been partly supported by the Czech Republic under
projects MSM6840770001 and AV0Z20760514. A part of this work was done while
Jakub~\v{S}\'{\i}stek was visiting at the University of Colorado Denver.
Another part of the work was performed during the stay of Jakub \v S\'\i stek
and Marta \v Cert\'\i kov\' a at Edinburgh Parallel Computing Center funded by
HPC-Europa project (RII3-CT-2003-506079). The authors would also like to thank
Bed\v{r}ich Soused\'{\i}k for his help in both typesetting and proofreading of
the manuscript.

% use this with bibtex - link to bibliography bddc.bib
%\bibliographystyle{elsart-num-sort}
%\bibliography{/home/sistja/denver/bddc/bibliography/bddc}

% included bbl file:

\begin{table}[p]
\begin{center}
\begin{tabular}
[c]{|c|c|c|c|c|}\hline
coarse problem & corners & corners+edges & corners+faces &
corners+edges+faces\\\hline\hline
iterations & 38 & 19 & 17 & 13\\\hline
cond. number est. & 117 & 15 & 65 & 7\\\hline\hline
factorization (sec) & 49 & 56 & 52 & 57\\\hline
pcg iter (sec) & 21 & 11 & 10 & 8\\\hline\hline
total (sec) & 85 & 85 & 80 & 83\\\hline
\end{tabular}
\newline\ \newline
\end{center}
\caption{Cube problem, $8$ subdomains, adding averages, initial set of
corners}
\label{tab:cube8}
\end{table}

\begin{table}[p]
\begin{center}
\begin{tabular}
[c]{|c|c|c|c|c|}\hline
coarse problem & corners & corners+edges & corners+faces &
corners+edges+faces\\\hline\hline
iterations & 42 & 16 & 24 & 11\\\hline
cond. number est. & 55 & 8 & 27 & 4\\\hline\hline
factorization (sec) & 2.2 & 3.2 & 2.8 & 4.0\\\hline
pcg iter (sec) & 7.5 & 3.7 & 5.1 & 3.6\\\hline\hline
total (sec) & 11.6 & 8.8 & 10.4 & 9.6\\\hline
\end{tabular}
\newline\ \newline
\end{center}
\caption{Cube problem, $64$ subdomains, adding averages, initial set of
corners}
\label{tab:cube64}
\end{table}

\begin{figure}[p]
\centering  \hbox{ \vbox{ \hsize=70mm \centering
\includegraphics[width=70mm]{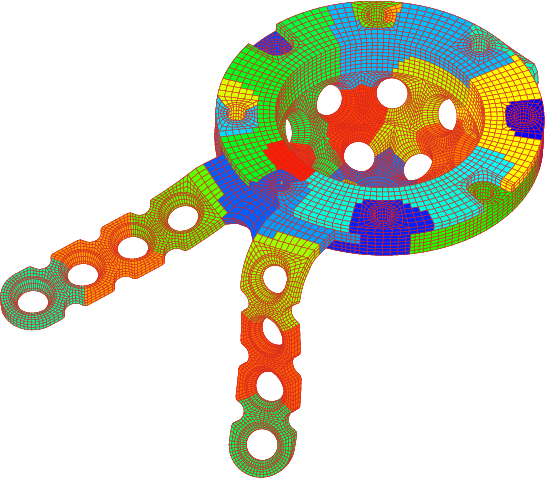}
} \vbox{ \hsize=70mm \centering
\includegraphics[width=70mm]{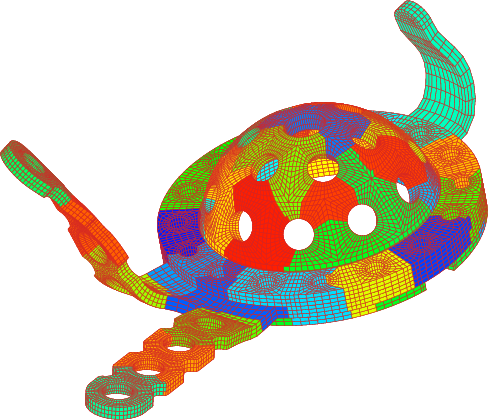}
}
}\caption{Hip joint replacement, division into 32 subdomains}
\label{fig:hip}
\end{figure}

\begin{table}[p]
\begin{center}
\begin{tabular}
[c]{|c|c|c|c|c|}\hline
coarse problem & corners & corners+edges & corners+faces &
corners+edges+faces\\\hline\hline
iterations & 181 & 171 & 69 & 62\\\hline
cond. number est. & 4,391 & 3,760 & 535 & 522\\\hline\hline
factorization (sec) & 100 & 70 & 86 & 80\\\hline
pcg iter (sec) & 241 & 216 & 94 & 87\\\hline\hline
total (sec) & 380 & 321 & 216 & 203\\\hline
\end{tabular}
\newline\ \newline
\end{center}
\caption{Hip joint replacement, $16$ subdomains, adding averages, 35 corners}
\label{tab:hip16basic}
\end{table}

\begin{table}[p]
\begin{center}
\begin{tabular}
[c]{|c|c|c|c|c|}\hline
coarse problem & corners & corners+edges & corners+faces &
corners+edges+faces\\\hline\hline
iterations & $>$500 & $>$500 & 137 & 70\\\hline
cond. number est. & n/a & n/a & n/a & n/a\\\hline\hline
factorization (sec) & 79 & 65 & 53 & 52\\\hline
pcg iter (sec) & $>$545 & $>$547 & 166 & 87\\\hline\hline
total (sec) & $>$651 & $>$638 & 236 & 161\\\hline
\end{tabular}
\newline\ \newline
\end{center}
\caption{Hip joint replacement, $32$ subdomains, adding averages, 57 corners}
\label{tab:hip32basic}
\end{table}

\begin{table}[ptb]
\begin{center}
\begin{tabular}
[c]{|c|c|c|c|c|}\hline
coarse problem & corners & corners+edges & corners+faces &
corners+edges+faces\\\hline\hline
iterations & 35 & 34 & 26 & 26\\\hline
cond. number est. & 96 & 96 & 65 & 65\\\hline\hline
factorization (sec) & 91 & 80 & 78 & 106\\\hline
pcg iter (sec) & 53 & 49 & 38 & 37\\\hline\hline
total (sec) & 183 & 166 & 153 & 181\\\hline
\end{tabular}
\newline\ \newline
\end{center}
\caption{Hip joint replacement, $16$ subdomains, adding averages, 335 corners}
\label{tab:hip16optimal}
\end{table}

\begin{table}[ptb]
\begin{center}
\begin{tabular}
[c]{|c|c|c|c|c|}\hline
coarse problem & corners & corners+edges & corners+faces &
corners+edges+faces\\\hline\hline
iterations & 35 & 32 & 30 & 27\\\hline
cond. number est. & 149 & 70 & 59 & 46\\\hline\hline
factorization (sec) & 60 & 57 & 59 & 62\\\hline
pcg iter (sec) & 49 & 40 & 37 & 34\\\hline\hline
total (sec) & 128 & 115 & 113 & 113\\\hline
\end{tabular}
\newline\ \newline
\end{center}
\caption{Hip joint replacement, $32$ subdomains, adding averages, 557 corners}
\label{tab:hip32optimal}
\end{table}

\end{document}